\documentclass[a4paper, 12pt]{amsart}

\usepackage{amsfonts}                   
\usepackage{latexsym}                   
\usepackage{amssymb}
\usepackage{mathptmx}  
\usepackage{graphics}
\usepackage[all]{xy}     
\newtheorem{teorema}{Theorem}
\newtheorem{Lemma}[teorema]{Lemma}
\newtheorem{propos}[teorema]{Proposition}
\newtheorem{corol}[teorema]{Corollary}
\newtheorem{ex}{Example}[section]
\newtheorem{rem}{Remark}[section]
\newtheorem{defin}[teorema]{Definition}

\def\defin{\par\ifdim\lastskip<\smallskipamount\removelastskip
  \smallskip\fi\noindent{\bf\ignorespaces
Definition\unskip:\enspace}\rm \ignorespaces}

\def\bit{\begin{itemize}}
\def\eit{\end{itemize}}
\def\be{\begin{equation}}
\def\ee{\end{equation}}
\def\beq{\begin{eqnarray}}
\def\eeq{\end{eqnarray}}

\def\ba{\begin{array}}
\def\ea{\end{array}}
\def\bt{\begin{teorema}}
\def\et{\end{teorema}}
\def\bp{\begin{propos}}
\def\ep{\end{propos}}
\def\bl{\begin{lemma}}
\def\el{\end{lemma}}
\def\bc{\begin{corol}}
\def\ec{\end{corol}}
\def\br{\begin{rem}\rm}
\def\er{\end{rem}}
\def\bex{\begin{ex}\rm}
\def\eex{\end{ex}}
\def\bd{\begin{defin}}
\def\ed{\end{defin}}
\def\demo{\par\noindent{\bf Proof.\ }}
\def\enddemo{\ $\Box$\par\vskip.6truecm}

     \def\nin{\noindent}

\def\ES{\varnothing}  \def\IN{\infty}  

\def\R{{\mathbb {R}}}    
\def\N{{\mathbb N}}     \def\d {\delta} \def\e{\varepsilon}
\def\C{{\mathbb C}}      
\def\Q{{\mathbb Q}}					
\def\Z{{\mathbb Z}}					\def\s{\sigma}

  \def\smi{\smallsetminus}

\def\sbs{\!\subset\!}

\def\oli{\overline}

\def\benu{\begin{enumerate}} \def\eenu{\end{enumerate}}
\def\beqn{\begin{eqnarray*}}  \def\eeqn{\end{eqnarray*}}

\def\beqn{\begin{eqnarray*}}  \def\eeqn{\end{eqnarray*}}

\begin{document}
\title{Cohomology and extension problems for semi $q$-coronae}

\author{Alberto Saracco$^{\star}$,\\ Giuseppe Tomassini$^{*,\star}$}
\thanks{$^{\star}$ Supported by the MURST project \lq\lq Geometric Properties of Real and
Complex Manifolds".\\
$^*$ Scuola Normale Superiore, Piazza dei Cavalieri, 7 - I-56126 Pisa, Italy.\\ \emph{e-mail: a.saracco@sns.it, g.tomassini@sns.it}}
\maketitle
\begin{abstract}
We prove some extension theorems for analytic objects, in particular sections of a coherent sheaf, defined in semi $q$-coronae of a complex space. Semi $q$-coronae are domains whose boundary is the union of a Levi flat part, a $q$-pseudoconvex part and a $q$-pseudoconcave part. Such results are obtained mainly using cohomological techniques.
\end{abstract}

\section{Introduction.}\nin
Let $X$ be a (connected and reduced) complex space. We
recall that $X$ is said to be {\it strongly} $q$-{\it pseudoconvex} in the sense of
Andreotti-Grauert~\cite{AG} if there exists a compact subset $K$ and a smooth function
$\varphi:X\to\R$, $\varphi\ge 0$, which is strongly $q$-plurisubharmonic on $X\smi
K$ and such that:  
\begin{itemize}
 \item[a)] $0=\min\limits_X\,\varphi<\min\limits_K\,\varphi$;
\item[b)] for every $c>\max\limits_K\,\varphi$ the subset
$$ 
B_c=\{x\in X:\varphi(x)<c\}  
$$  
is relatively compact in $X$.  
\end{itemize}
If $K=\ES$, $X$ is said to be $q$-{\it complete}. We remark that, for a space, being $1$-complete is equivalent to being Stein.

Replacing the condition b) by 
\bit
\item[b')] for every $0<\varepsilon<\min\limits_K\,\varphi$
and $c>\max\limits_K\,\varphi$ the subset
$$ 
B_{\e,c}=\{x\in X:\e<\varphi(x)<c\}  
$$
is relatively compact in $X$,
\eit
we obtain the notion of $q$-{\it corona} (see~\cite{AG},~\cite{AT}).

A $q$-corona is said to be {\it
complete} whenever $K=\ES$.     

The extension problem for analytic objects defined on $q$-coronae was studied by many authors (see e.g. \cite{FG}, \cite{Se}, \cite{Si}, \cite{SiT}, \cite{T70}). In this paper we deal with the larger class of the
semi $q$-coronae which are defined as follows. Consider a strongly $q$-pseudoconvex space (or,
more generally, a $q$-corona) $X$, and a smooth function $\varphi:X\to\R$ displaying the
$q$-pseudoconvexity of $X$. Let $B_{\varepsilon,c}\sbs X$ and let $h:X\to\R$ be a
pluriharmonic function (i.e.\ locally the real part of a holomorphic function)
such that $K\cap\{h=0\}=\ES$. A connected component of $B_{\varepsilon,c}\smi\{h=0\}$ is, by definition, a {\it
semi} $q$-{\it corona}. 

Another type of semi $q$-corona is obtained by replacing the zero set of $h$ with the
intersection of $X$ with a Levi flat hypersurface. More precisely, consider a closed
strongly $q$-pseudoconvex subspace $X$ of an open subset of $\C^n$ and the $q$-corona 
$C=B_{\varepsilon,c}=B_c\smi{\overline B}_\varepsilon$. Let $H$ be a Levi flat
hypersurface of a neighbourhood $U$ of ${\overline B}_c$ such that $H\cap K=\ES$. The
connected components $C_m$ of $C\smi H$ are called semi $q$-coronae.

In both cases the semi $q$-coronae are differences $A_c\smi{\overline A}_\varepsilon$
where $A_c$, $A_\varepsilon$ are strongly $q$-pseudoconvex spaces. Indeed, the function
$\psi=-\log h^2$ (respectively $\psi=-\log \delta_H(z)$, where $\delta_H(z)$ is the
distance of $z$ from $H$) is plurisubharmonic in $W\smi\{h=0\}$ (respectively 
$W\smi H$) where $W$ is a  neighbourhood of $B_c\cap \{h=0\}$ (respectively $B_c\cap H$). Let $\chi:\R\to\R$ be an
increasing convex function such that $\chi\circ\varphi>\psi$ on a neighbourhood of
$B_c\smi W$. The function $\Phi=\sup\,(\chi\circ\varphi,\psi)+\varphi$ is an exhaustion function 
for $B_c\smi \{h=0\}$ (for $B_c\smi H$) and it is strongly
$q$-plurisubharmonic in $B_c\smi(\{h=0\}\cup K)$ (in 
$B_c\smi H\cup K)$.

The interest for
domains whose boundary contains a ``Levi flat part'' originated from an extension theorem for
CR-functions proved in \cite{LT} (see also \cite{La}, \cite{LaP}, \cite{St}).   

Using cohomological techniques developped in \cite{AG}, \cite{AT}, \cite{BS}, \cite{C} we prove that, under
appropriate regularity conditions, holomorphic functions defined on a complete semi
$1$-corona \lq\lq fill in the holes\rq\rq\ (Corollaries~\ref{cD} and~\ref{oE}). Meanwhile we also
obtain more general extension theorems for sections of coherent  sheaves
(Theorems~\ref{cC} and~\ref{oCw}). As an application, we finally obtain an
extension theorem for divisors (Theorems~\ref{divis} and~\ref{divis2}) and for analytic sets
of codimension one (Theorem~\ref{ansets}).

We remark that this approach fails in the case when the objects to be extended are not
sections of a sheaf defined on the whole $B_c$. In particular, this applies for analytic
sets of higher codimension. This is closely related with the general, definitely more
difficult, problem of extending analytic objects assigned on some semi $q$-corona when
the subsets $B_c$ are not relatively compact in $X$ i.e.\ when $X$ is a genuine $q$-corona. It is worth noticing that a similar extension theorem for complex submanifold of higher codimension has been recently obtained in~\cite{DS} by different methods based on Harvey-Lawson's theorem~\cite{HL}.\vspace{0,3cm}

We wish to thank Mauro Nacinovich and Viorel V\^aj\^aitu for their kind help and suggestions.

\section{Cohomology and extension of sections.}\nin
\subsection{Closed $q$-coronae}

Let $X$ be a strictly $q$-pseudoconvex space (respectively $X\subset\C^n$ be a strictly $q$-pseudoconvex open set) and $H=\left\{h=0\right\}$ (respectively $H$ Levi-flat), and $C=B_{\varepsilon,c}=B_c\smi{\overline B}_\varepsilon$ a $q$-corona.

We can suppose that $B_c \smi H$ has two connected components, $B_+$ and $B_-$, and define $C_+=B_+\cap C$, $C_-=B_-\cap C$.

If $\mathcal F\in{\rm Coh}(B_c)$, we define $p(\mathcal F)=\inf\limits_{x\in
B_c}\,{\rm depth}({\mathcal F}_x)$, the depth of $\mathcal F$ on $B_c$. If $\mathcal F=\mathcal O$, the structure sheaf of $X$, we define $p(B_c)=p(\mathcal O)$.

\bt\label{Ac}
Let $\mathcal F\in {\rm Coh}(B_c)$. Then the image of the homomorphism 
$$
H^r(\oli{B}_+,\mathcal F)\oplus H^r(\oli{C},\mathcal F)\longrightarrow H^r(\oli{C}_+,\mathcal F)
$$
(all closures are taken in $B_c$), defined by $(\xi\oplus\eta)\mapsto\xi_{|\oli{C}_+}-\eta_{|\oli{C}_+}$ has finite codimension provided that $q-1\le
r\le p(\mathcal F)-q-2$ . 
\et
\demo
Consider the Mayer-Vietoris sequence applied to the closed sets $\oli{B}_+$ and $\oli{C}$
\beq\label{suc1}
\cdots &\to& H^r(\oli{B}_+\cup \oli{C},\mathcal F)\to H^r(\oli{B}_+,\mathcal F)\oplus H^r(\oli{C},\mathcal
F)\stackrel{\d}{\to}\\ &\stackrel{\d}{\to}& H^r(\oli{C}_+,\mathcal F)\to H^{r+1}(\oli{B}_+\cup \oli{C},\mathcal
F)\to\cdots\nonumber 
\eeq
$\d(a\oplus b)=a_{|\oli{C}_+}-b_{|\oli{C}_+}$.
$\oli B_+\cup \oli C=B_c\smi U$ where $U=B_-\cap B_\varepsilon$. $U$ is $q$-complete, so the groups of compact support cohomology $H^{r}_c(U,\mathcal F)$ are zero for $q\leq r\leq p(\mathcal{F})-q$.

From the exact sequence of compact support cohomology
\beq
\cdots &\to& H^r_c(U,\mathcal F)\to H^r(B_c,\mathcal F)\to\\ 
&\to& H^r(B_c\smi U,\mathcal F)\to H^{r+1}_c(U,\mathcal F)\to\cdots\nonumber  
\eeq
it follows that
\begin{equation}\label{isomBc-U}
H^r(B_c,\mathcal F)\stackrel{\sim}{\rightarrow} H^r(B_c\smi U,\mathcal F),
\end{equation}
for $q\leq r \leq p(\mathcal{F})-q-1$.

Since $B_c$ is $q$-pseudoconvex,
$$\dim_\C\,H^r(B_c,\mathcal F)<\IN
$$
for $q\le r$ \cite[Th\'eor\`eme 11]{AG}, and so
$$
\dim_\C\,H^r( B_c\smi U,\mathcal F)<\IN
$$
for $q\le r\le p(\mathcal F)-q-1$.

From (\ref{suc1}) we see that $\dim_\C H^r( B_c\smi U,\mathcal F)=\dim_\C H^r(\oli B_+\cup\oli C ,\mathcal F)$ is greater than or equal to the codimension of the homomorphism $\delta$.
\enddemo
\bc\label{cB}
Under the same assumption of Theorem~\ref{Ac}, if $K\cap H=\ES$,
$$
\dim_\C\,H^r(\oli C_+,\mathcal F)<\IN
$$
for $q\le r\le p(\mathcal F)-q-2$.
\ec
\demo
Since $K\cap H=\ES$, $\oli B_+$ is a $q$-pseudoconvex space, and by virtue of \cite[Th\'eor\`eme 11]{AG} we have
$$
\dim_\C\,H^r(\oli B_+,\mathcal F)<\IN
$$ 
for $r\ge q$. On the other hand, $\oli C$ is
a $q$-corona, thus we obtain
$$
\dim_\C\,H^r(\oli C,\mathcal F)<\IN
$$
for $q\le r\le p(\mathcal F)-q-1$ in view of \cite[Theorem 3]{AT}. By Theorem~\ref{Ac} we then get that for $q\le r\le p(\mathcal
F)-q-1$ the vector space $H^r(\oli B_+,\mathcal F)\oplus H^r(\oli C,\mathcal F)$ has finite dimension and for $q-1\leq r\leq p(\mathcal F)-q-2$ its image in $H^r(\oli C_+,\mathcal F)$ has finite codimension. Thus 
$H^r(\oli C_+,\mathcal F)$ has finite dimension for $q\le r\le p(\mathcal F)-q-2$.
\enddemo
\bt\label{cC}
If $\oli B_+$ is a $q$-complete space, then
$$
H^r(\oli C,\mathcal F)\stackrel{\sim}{\rightarrow} H^r(\oli C_+,\mathcal F)
$$
for $q\le r\le p(\mathcal F)-q-2$ and the homomorphism
\begin{equation}\label{eqA}
H^{q-1}(\oli B_+,\mathcal F)\oplus H^{q-1}(\oli C,\mathcal F)\longrightarrow H^{q-1}(\oli C_+,\mathcal F)
\end{equation}
is surjective for $p(\mathcal F)\geq2q+1$.

If $\oli B_+$ is a $1$-complete space and $p(\mathcal F)\ge 3$, the homomorphism
$$
H^0(\oli B_+,\mathcal F)\longrightarrow H^0(\oli C_+,\mathcal F)
$$
is surjective. 
\et
\demo
Since by hypothesis $\oli B_+$ is a $q$-complete space, $H^r(\oli B_+,\mathcal F)=\{0\}$ for
$q\le r$ \cite[Th\'eor\`eme 5]{AG}. From (\ref{isomBc-U}) it follows that $H^r(\oli B_+\cup \oli C_+,\mathcal F)=\{0\}$ for 
$q\le r\le p(\mathcal F)-q-1$. Thus, the Mayer-Vietoris sequence (\ref{suc1}) implies that $H^r(\oli C,\mathcal F)\stackrel{\sim}{\rightarrow} H^r(\oli C_+,\mathcal F)$ for $q \le r \le p(\mathcal F)-q-2$ and that the homomorphism (\ref{eqA}) is surjective if $p(\mathcal F)\geq 2q+1$.

In particular, if $q=1$ and $p(\mathcal F)\ge 3$ the homomorphism$$
H^0(\oli B_+,\mathcal F)\oplus H^0(\oli C,\mathcal F)\longrightarrow H^0(\oli C_+,\mathcal F)
$$
is surjective i.e.\ every section $\s\in H^0(\oli C_+,\mathcal F)$ is a difference $\s_1-\s_2$ of two sections $\s_1\in H^0(\oli B_+,\mathcal F)$, $\s_2\in H^0(\oli C,\mathcal F)$. Since $B_\varepsilon$ is Stein, the cohomology group with compact supports $H^1_k(B_\varepsilon,\mathcal F)$
is zero, and so the Mayer-Vietoris compact support cohomology sequence implies that the restriction homomorphism 
$$ 
H^0(\oli B_c,\mathcal F)\longrightarrow H^0(\oli B_c\smi B_\varepsilon,\mathcal
F)=H^0(\oli C,\mathcal F)
$$
is surjective, hence $\s_2\in H^0(\oli C,\mathcal F)$ is restriction of $\widetilde{\s}_2\in H^0(B_c,\mathcal F)$. So $\s$ is restriction to $\oli C_+$ of $(\s_1-\widetilde{\s}_{2|\oli B_+})\in H^0(\oli B_+,\mathcal F)$, and the restriction homomorphism is surjective. 
\enddemo
\bc\label{cD}
Let $\oli B_+$ be a $1$-complete space and $p(B_c)\ge 3$. Then every holomorphic function on $\oli C_+$ extends holomorphically on $\oli B_+$.
\ec

\subsection{Open $q$-coronae}

Most of the Theorems and Corollaries of the previous section still hold in the open case and their proofs are very similar. First we give the proof of the extension results using directly Theorem~\ref{cC}. We have to assume that $H$ is the zero set of a pluriharmonic function $h$ and we define $B_c$, $C$, $B_+$, $B_-$, $C_+$ and $C_-$ as we did before.

Let us suppose $B_+$ is $1$-complete and $p(\mathcal F)\geq 3$. Let $s\in H^0(C_+,\mathcal F)$. For all $\epsilon>0$, we consider the closed semi $1$-corona
$$
\oli C_\epsilon=\oli{B_{\varepsilon+\epsilon,c}\cap\{h>\epsilon\}}\subset C_+
$$
Let $\s_\epsilon=s_{|\oli C_\epsilon}$. By Theorem~\ref{cC} (applied to $
\oli C_\epsilon$, $H_\epsilon=\{h=\epsilon\}$), we obtain that $\s_\epsilon$ extends to a section $\widetilde\s_\epsilon\in H^0(\oli B_\epsilon,\mathcal F)$, where $\oli B_\epsilon=\oli{B_+\cap\{h>\epsilon\}}$. Since $B_+=\cup_\epsilon \oli B_\epsilon$, if for all $\epsilon_2>\epsilon_1>0$,
\beq\label{*}
\widetilde\s_{\epsilon_1|_{\oli B_{\epsilon_2}}}=\widetilde\s_{\epsilon_2}
\eeq
the sections $\widetilde\s_{\epsilon}$ can be glued toghether to a section $\s\in H^0(B_+,\mathcal F)$ extending $s$.

Let $\epsilon_1,\epsilon_2$, $\epsilon_2>\epsilon_1>0$, be fixed. We have to show that~(\ref{*}) holds. By definition,
$$
\left(\widetilde\s_{\epsilon_1|_{\oli B_{\epsilon_2}}}-\widetilde\s_{\epsilon_2}\right)_{|\oli C_{\epsilon_2}}=s-s=0.
$$
Thus, the support of $\widetilde\s_{\epsilon_1|_{\oli B_{\epsilon_2}}}-\widetilde\s_{\epsilon_2}$, $S$, is an analytic set contained in $\oli B_{\epsilon_2}\smi C_{\epsilon_2}$. Let us consider the family $(\phi_\lambda=\lambda(\varphi-\epsilon_2)+(1-\lambda)(h-\epsilon_2))_{\lambda\in[0,1]}$ of strictly plurisubharmonic functions. Let $\oli \lambda$ be the smallest value of $\lambda$ for which $\{\phi_\lambda=0\}\cap S\neq\ES$. Then $\{\phi_{\oli \lambda}<0\}\cap B_+\subset B_+$ is a Stein domain in which the analytic set $S$ intersects the boundary; so the maximum principle for plurisubharmonic functions and the strict plurisubharmonicity of $\phi_{\oli \lambda}$ toghether imply that $\{\phi_{\oli \lambda}=0\}\cap S$ is a set of isolated points in $S$. By repeating the argument, we show that $S$ has no components of positive dimension. Hence $\widetilde\s_{\epsilon_1|_{\oli B_{\epsilon_2}}}-\widetilde\s_{\epsilon_2}$ is zero outside a set of isolated points. Since $p(\mathcal F)\geq3$, the only section of $\mathcal F$ with compact support is the zero-section \cite[Th\'eor\`eme 3.6 (a), p.\ 46]{BS}, and so $\widetilde\s_{\epsilon_1|_{\oli B_{\epsilon_2}}}-\widetilde\s_{\epsilon_2}$ is zero.

Hence, there exists a section $\s\in H^0(B_+,\mathcal F)$ such that $\s_{|C_+}=s$. Thus we have proved the following

\bt\label{oCw}
If a $B_+$ is $1$-complete space, $\mathcal F$ a coherent sheaf on $B_+$ with $p(\mathcal F)\ge 3$, the homomorphism
$$
H^0(B_+,\mathcal F)\longrightarrow H^0(C_+,\mathcal F)
$$
is surjective.
\et

In particular,
\bc\label{oE}
If $B_+$ is a $1$-complete space and $p(B_c)\ge 3$, every holomorphic function on $C_+$
can be holomorphically extended on $B_+$.
\ec


\bt\label{oA}
Let $\emph{Sing}(B_c)=\ES$. Let $\mathcal F\in {\rm Coh}(B_c)$. Then the image of the homomorphism 
$$
H^r(B_+,\mathcal F)\oplus H^r(C,\mathcal F)\longrightarrow H^r(C_+,\mathcal F)
$$
defined by $(\xi,\eta)\mapsto\xi_{|C_+}-\eta_{|C_+}$ has finite codimension for $q-1\le
r\le p(\mathcal F)-q-2$. For $q=1$ the thesis holds true also dropping the assumption $\emph{Sing}(B_c)=\ES$.
\et

\demo
Consider the Mayer-Vietoris sequence applied to the open sets $B_+$ and $C$
\beq\label{1open}
\cdots &\to& H^r(B_+\cup C,\mathcal F)\to H^r(B_+,\mathcal F)\oplus H^r(C,\mathcal
F)\stackrel{\d}{\to}\\ &\stackrel{\d}{\to}& H^r(C_+,\mathcal F)\to H^{r+1}(B_+\cup C,\mathcal
F)\to\cdots,\nonumber 
\eeq
$\d(a\oplus b)=a_{|C_+}-b_{|C_+}$. $B_+\cup C=B_c\smi K_0$ where $K_0=\oli{B}_-\cap\oli{B}_\varepsilon$. $K_0$ has a $q$-complete neighbourhoods system and so the local cohomology groups $H^r_{K_0}(B_c,\mathcal F)$ are zero for $q\leq r\le p(\mathcal F)-q$ \cite{C} (in the general case for $q=1$, see \cite[Lemme 2.3, p.\ 29]{BS}). 

Then, from the local cohomology exact sequence
\beq
\cdots &\to& H^r_{K_0}(B_c,\mathcal F)\to H^r(B_c,\mathcal F)\to\\ 
&\to& H^r(B_c\smi K_0,\mathcal F)\to H^{r+1}_{K_0}(B_c,\mathcal F)\to\cdots\nonumber  
\eeq
follows that 
\begin{equation}\label{eq3}
H^r(B_c,\mathcal F)\stackrel{\sim}{\rightarrow} H^r(B_c\smi K_0,\mathcal F),
\end{equation}
for $q\le r\le p(\mathcal F)-q-1$.

Since $B_c$ is $q$-pseudoconvex,
$$
\dim_\C\,H^r(C,\mathcal F)<\IN
$$
for $q\le r$ \cite[Th\'eor\`eme 11]{AG}, and so
$$
\dim_\C\,H^r(B_c\smi K_0,\mathcal F)<\IN
$$
for $q\le r\le p(\mathcal F)-q-1$.

From (\ref{1open}) we see that $\dim_\C H^r(B_c\smi K_0,\mathcal F)=\dim_\C H^r(B_+\cup C,\mathcal F)$ is greater than or equal to the codimension of the homomorphism $\d$.
\enddemo

\bc\label{oB}
Under the same assumption of Theorem~\ref{oA}, if $K\cap H=\ES$,
$$
\dim_\C\,H^r(C_+,\mathcal F)<\IN
$$
for $q\le r\le p(\mathcal F)-q-2$.
\ec
\demo
The proof is similar to that of Corollary~\ref{cB}.
\enddemo

\bt\label{oC}
Suppose that $\emph{Sing}(B_c)=\ES$ and $B_+$ is a $q$-complete space, then
$$
H^r(C,\mathcal F)\stackrel{\sim}{\rightarrow} H^r(C_+,\mathcal F)
$$
for $q\le r\le p(\mathcal F)-q-2$ and the homomorphism
\beq\label{Aopen}
H^{q-1}(B_+,\mathcal F)\oplus H^{q-1}(C,\mathcal F)\longrightarrow H^{q-1}(C_+,\mathcal F)
\eeq
is surjective if $p(\mathcal F)\geq 2q+1$. If $q=1$, both results hold true for an arbitrary complex space $B_c$.

\et
\demo
The proof is similar to that of Theorem~\ref{cC}.
%
\enddemo

\subsection{Corollaries of the extension theorems.}
From now on, unless otherwise stated, by $B$, $B_+$, $B_\varepsilon$, $C$ and  $C_+$ we denote both the open sets and their closures, and we suppose that $H=\{h=0\}$, $h$ pluriharmonic.

\subsubsection{}Let $f\in H^0(C_+,\mathcal O^*)$. In the hypothesis of Corollaries~\ref{cD} and~\ref{oE}, both $f$ and $1/f$ extend holomorphically on $B_+$ Hence:
\bc\label{O*sur}
If $B_+$ is a $1$-complete space and $p(B_c)\ge 3$, the restriction homomorphism
$$
H^0(B_+,\mathcal O^*)\longrightarrow H^0(C_+,\mathcal O^*)
$$
is surjective.
\ec

\subsubsection{}In Theorems~\ref{cC} and~\ref{oCw} we have estabilished the isomorphism
$$
H^r(C,\mathcal F)\stackrel{\sim}{\rightarrow} H^r(C_+,\mathcal F).
$$
In some special cases this leads to vanishing-cohomology theorems for $C_+$. An example is provided by a $q$-corona $C$ which is contained in an affine variety. In such a situation, we have that $H^r(C,\mathcal F)=\left\{0\right\}$, for $q\leq r\leq p(\mathcal F)-q-2$ \cite{AT}, and consequently $H^r(C_+,\mathcal F)=\left\{0\right\}$ in the same range of $r$.

\subsubsection{}Let $X$ be a Stein space. Let $H=\left\{h=0\right\}\subset X$ be the zero set of a pluriharmonic function, and let $S$ be a real hypersurface of $X$ with boundary, such that $S\cap H=b S=b A$, where $A$ is an open set in $H$. Let $D\subset X$ be the relatively compact domain bounded by $S\cup A$. In \cite{LT} it is proved that, for $X=\C^n$, $CR$-functions on $S$ extend holomorphically to $D$. As a corollary of the previous theorems, we can obtain a similar result for section of a coherent sheaf on an arbitrary Stein space $X$.

Let us consider the connected component $Y$ of $X\smi H$ containing $D$, the closure $\oli D$ of $D$ in $Y$, and let be $F=Y\smi D$ and $S_Y=S\cap Y$. For every coherent sheaf $\mathcal F$ on $X$, with $p(\mathcal F)\geq3$ we have the Mayer-Vietoris exact sequence
$$
\cdots\ \to\ H^0(\oli D,\mathcal F)\oplus H^0(F,\mathcal F)\ \to\ H^0(S_Y,\mathcal F)\ \to\ H^1(Y, \mathcal F) \ \to\ \cdots
$$
Since $Y$ is Stein, $H^1(Y, \mathcal F)$ is zero, and every section $\s$ on $S_Y$ is a difference $s_1-s_2$, where $s_1\in H^0(\oli D,\mathcal F)$ and $s_2\in H^0(F,\mathcal F)$. By choosing an $\varepsilon$ big enough so that $S$ is contained in the ball $B_\e(x_0)$ of radius $\e$ of $X$ centered in $x_0$, we can apply Theorem~\ref{oCw} to the semi $1$-corona $C_+=Y\smi(B_\e\cap Y)$, to extend $s_{2|_{C_+}}$ to a section $\tilde s_2$ defined on $Y$. In order to conclude that $s_1-\tilde s_{2|_{\oli D}}$ extends the section $\s$, we have to prove that $s_{2|_F}-\tilde s_{2|_F}=0$.

As before, we consider the set $\Sigma=\left\{s_{2|_F}-\tilde s_{2|_F}\neq0\right\}\subset B_\e\cap Y$ and conclude that $\Sigma$ is a set of isolated points. Since $p(\mathcal F)\geq3$, $\mathcal F$ has no non zero section with compact support \cite[Th\'eor\`eme 3.6 (a), p.\ 46]{BS}. Thus $\Sigma=\ES$ and we have obtained the following:

\bc\label{Lupac1}
Let $X$ be a Stein space. Let $H=\left\{h=0\right\}\subset X$ be the zero set of a pluriharmonic function, and $S$ be a real hypersurface of $X$ with boundary, such that $S\cap H=b S=b A$, where $A$ is an open set in $H$. Let $D\subset X$ be the relatively compact domain bounded by $S\cup A$ and $\mathcal F$ be a coherent sheaf with $p(\mathcal F)\geq3$. All sections of $\mathcal F$ on $S$ extend (uniquely) to $D$.
\ec

We can go further:

\bc\label{Lupac} Let $X$ be a Stein manifold, $\mathcal F$ a coherent sheaf on $X$ such that $p(\mathcal F)\geq3$ and $D$ be a bounded domain and $K$ a compact subset of $b D$ such that $b D\smi K$ is smooth. Assume that $K$ is $\mathcal O(D)$-convex, i.e.
$$ K=\left\{z\in\oli D\ :\ |f(z)|\leq\max_K |f|\right\}. $$
Then every section of $\mathcal F$ on $b D\smi K$ extends to $D$.
\ec
\demo We recall that since $U$ is an open subset of a Stein manifold there exists
an envelope of holomorphy $\widetilde U$ of $U$ (cfr. \cite{DG}) $\widetilde U$ is a Stein
domain $\pi_U:\widetilde U\to X$ over $X$ and there exists and open embedding $j:U\to\widetilde U$ such that $\pi_U\circ j=id_U$ and $J^*:\mathcal O(\widetilde U)\to\mathcal O(U)$ is an isomorphism. In particular $\pi_U^*\mathcal F$ is a coherent sheaf with
the same depth as $\mathcal F$, which extends ${\mathcal F}_{|U}$.

Let us fix an arbitrary point $x\in D$. We need to show that any given section $\s\in H^0(b D\smi K,\mathcal F)$ extends to a neighbourhood of $x$. Since $x\not\in K=\widehat{K}$, there exists an holomorphic function $f$, defined on a neighbourhood $U$ of $\oli D$, such that $|f(x)|>\max_K |f(z)|$.

Then $\s$ extends to a section $\widetilde\s\in H^0(\pi^{-1}(D\smi K),\mathcal F)$. Let $\widetilde f$ be the
holomorphic extension of $f$ to $\widetilde U$. The hypersurface  
$$
H=\left\{z\in\widetilde U:\vert \widetilde
f(z)\vert=\max\limits_K\vert\widetilde f\vert\right\}
$$
is the zero-set of a pluriharmonic function and, by construction, $$x\in \widetilde D_+=\left\{z\in\widetilde U:\vert \widetilde
f(z)\vert>\max\limits_K\vert\widetilde f\vert\right\}.$$
Now we are in the situation of Corollary~\ref{Lupac1} so $\widetilde\s$ extends to a section on $\widetilde D_+$. Since 
$x\in\widetilde D_+$, this ends the proof.
\enddemo


\section{Extension of divisors and analytic sets of codimension one.}
First of all, we give an example in dimension $n=2$ of a regular complex curve of $C_+$ which does not extend on $B_+$. Hence, not every divisor on $C_+$ extends to a divisor on $B_+$.\vspace{0.25cm}\\

\nin\textbf{Example}. Using the same notation as before, let $B_c$ be the ball $\left\{|z_1|^2+|z_2|^2<c\right\}$, $c>2$, in $\C^2$, and $H$ be the hyperplane $\left\{x_2=0\right\}$ ($z_j=x_j+iy_j$). Let $2<\e<c$, $C=B_c\smi\oli B_\e$, $B_+=B_c\cap\left\{x_2>0\right\}$, $C_+=C\cap\left\{x_2>0\right\}$.

Consider the connected irreducible analytic set of codimension one
$$A=\{(z_1,z_2)\in B_+\ :\ z_1z_2=1\}$$
and its restriction $A_C$ to $C_+$. If $A_C$ has two connected components, $A_1$ and $A_2$, if we try to extend $A_1$ (analytic set of codimension one on $C_+$) to $B_+$, its restriction to $C_+$ will contain also $A_2$. So $A_1$ is an analytic set of codimension one on $C_+$ that does not extend on $B_+$.

So, let us prove that $A_C$ has indeed two connected components. A point of $A$ (of $A_C$) can be written as $z_1=\rho e^{i\theta}$, $z_2=\frac{1}{\rho} e^{-i\theta}$, with $\rho\in\R^+$ and $\theta\in\left(-\frac\pi2,\frac\pi2\right)$. Hence, points in $A_C$ satisfy
$$
2<\varepsilon<\rho^2+\frac1{\rho^2}<c\ \Rightarrow\ 2<\sqrt{\varepsilon+2}<\rho+\frac1\rho<\sqrt{c+2} .
$$
Since $f(\rho)=\rho+1/\rho$ is monotone decreasing up to $\rho=1$ (where $f(1)=2$), and then monotone increasing, there exist $a$ and $b$ such that the inequalities are satisfied when $a<\rho<b<1$, or when $1<1/b<\rho<1/a$. $A_C$ is thus the union of the two disjoint open sets
$$
\xymatrix{A_1=\left\{ \left(\rho e^{i\theta},\frac1\rho e^{-i\theta}\right)\in \C^2\ \Big|\ a<\rho<b,\  -\frac\pi2<\theta<\frac\pi2\right\};\\
A_2=\left\{ \left(\rho e^{i\theta},\frac1\rho e^{-i\theta}\right)\in \C^2\ \Big|\ a<\frac1\rho<b,\ -\frac\pi2<\theta<\frac\pi2\right\}.}$$\vspace{0.25cm}

The aim of this section is to prove an extension theorem for divisors, i.e.\ to prove that, under certain hypothesis, the homomorphism
\beq\label{DivSurg}
H^0(B_+,\mathcal D)\to H^0(C_+,\mathcal D)
\eeq
is surjective.

In order to get this result, we observe that from the exact sequence
\beq\label{esattaD}
0\to\mathcal O^*\to\mathcal M^*\to\mathcal D\to 0
\eeq
we get the commutative diagram (horizontal lines are exact)
$$
\xymatrix{H^0(B_+,\mathcal M^*)\ar[r]\ar[d]_{\alpha} & H^0(B_+,\mathcal D)\ar[r]\ar[d]_{\beta} & H^1(B_+,\mathcal O^*)\ar[r]\ar[d]_{\gamma} & H^1(B_+,\mathcal M^*)\ar[d]_{\delta} \\ H^0(C_+,\mathcal M^*)\ar[r] & H^0(C_+,\mathcal D)\ar[r] & H^1(C_+,\mathcal O^*)\ar[r] & H^1(C_+,\mathcal M^*)}
$$
Thus, in view of the \lq\lq five lemma\rq\rq, in order to conclude that $\beta$ is surjective it is sufficient to show that $\alpha$ and $\gamma$ are surjective, and $\delta$ is injective.

\begin{Lemma}\label{alphaS}
If $\emph{Sing}(B_+)=\ES$, $B_c$ is $1$-complete and $p(B_c)\geq 3$, then $\alpha$ is surjective.
\end{Lemma}
\demo
Let $f$ be a meromorphic invertible function on $C_+$. Since $C_+$ is an open set of the Stein manifold $B_+$, $f=f_1 f_2^{-1}$, $f_1,f_2\in H^0(C_+,\mathcal O)$. By Corollary~\ref{cD} (\ref{oE}), $f_1$ and $f_2$ extend to holomorphic functions on $B_+$ and consequently $f$ extends on $B_+$ as well.
\enddemo

\begin{Lemma}\label{gammaS}
Assume that the restriction $H^2(B_+,\Z)\to H^2(C_+,\Z)$ is surjective. If $B_c$ is $1$-complete and $p(B_c)\geq4$, then $\gamma$ is surjective.
\end{Lemma}
We remark that if $H^2(C_+,\Z)=\{0\}$ the first condition is satisfied.\vspace{0.5cm}
\demo
From the exact sequence
\beq
0\to\Z\to\mathcal O\to\mathcal O^*\to 0
\eeq
we get the commutative diagram (horizontal lines are exact)
$$
\xymatrix{H^1(B_+,\mathcal O)\ar[r]\ar[d]_{f_2} & H^1(B_+,\mathcal O^*)\ar[r]\ar[d]_{\gamma} & H^2(B_+,\Z)\ar[r]\ar[d]_{f_4} & H^2(B_+,\mathcal O)\ar[d]_{f_5} \\ H^1(C_+,\mathcal O)\ar[r] & H^1(C_+,\mathcal O^*)\ar[r] & H^2(C_+,\Z)\ar[r] & H^2(C_+,\mathcal O)}
$$
where $H^1(B_+,\mathcal O)=H^2(B_+,\mathcal O)=\left\{0\right\}$ because $B_+$ is Stein, and $f_4$ is surjective by hypothesis. Thus in order to prove that $\gamma$ is surjective by the \lq\lq five lemma\rq\rq\ it is sufficient to show that $f_2$ is surjective, i.e.\ that $H^1(C_+,\mathcal O)=\left\{0\right\}$.

Since $p(B_c)\geq4$, by Theorem~\ref{cC} (\ref{oC}) it follows that
\beq\label{C=C+}
H^1(C,\mathcal O)\stackrel{\sim}{\longrightarrow} H^1(C_+,\mathcal O).
\eeq
Consider the local, respectively compact support, cohomology exact sequence
$$
\xymatrix{H^1_{\oli B_\varepsilon}(B_c,\mathcal O)\ar[r] & H^1(B_c,\mathcal O)\ar[r] & H^1(C,\mathcal O)\ar[r] & H^2_{\oli B_\varepsilon}(B_c,\mathcal O) \\ H^1_k(B_\varepsilon,\mathcal O)\ar[r] & H^1(B_c,\mathcal O)\ar[r] & H^1(C,\mathcal O)\ar[r] & H^2_k(B_\varepsilon,\mathcal O)}
$$
Since $B_c$ is Stein, $H^1(B_c,\mathcal O)=\left\{0\right\}$ and $H^r_k(B_\e,\mathcal O)=H^r_{\oli B_\e}(B_c,\mathcal O)=\left\{0\right\}$ for $1\leq r\leq p(B_\e)-1$~\cite{C}. In particular, since $p(B_\e)\geq p(B_c)\geq4$, it follows that
\beq\label{Bc=C}
\{0\}=H^1(B_c,\mathcal O)\stackrel{\sim}{\longrightarrow} H^1(C,\mathcal O).
\eeq
(\ref{C=C+}) and (\ref{Bc=C}) give
$$
\{0\}=H^1(B_c,\mathcal O)\stackrel{\sim}{\longrightarrow}H^1(C,\mathcal O)\stackrel{\sim}{\longrightarrow} H^1(C_+,\mathcal O).
$$and this proves the lemma.
\enddemo

In the case $H^2(C_+,\Z)=\{0\}$ we remark that from the proof of Lemma~\ref{gammaS} it follows that the sequence 
$$\{0\}\longrightarrow H^1(C_+,\mathcal{O}^*)\longrightarrow\{0\}$$
is exact, that is $H^1(C_+,\mathcal{O}^*)=\{0\}$. Hence, the commutative diagram relative to (\ref{esattaD}) becomes (horizontal lines are exact)
\beq\label{last}
\xymatrix{H^0(B_+,\mathcal M^*)\ar[r]\ar[d]_{\alpha} & H^0(B_+,\mathcal D)\ar[r]\ar[d]_{\beta} & H^1(B_+,\mathcal O^*)\ar[d]_{\gamma} \\ H^0(C_+,\mathcal M^*)\ar[r] & H^0(C_+,\mathcal D)\ar[r] & \{0\}}
\eeq
and it is then easy to see that a divisor on $C_+$ can be extended to a divisor on $B_+$.

Thus we have proved the following:

\bt\label{divis}
Let $B_c$ be $1$-complete, $p(B_c)\geq4$, and $C_+$ satisfy the topological condition $H^2(C_+,\Z)=\{0\}$. Then, if $\emph{Sing}(B_+)=\ES$, all divisors on $C_+$ extend (uniquely) to divisors on $B_+$.
\et

\bc\label{corxi}
Let $B_c$ be $1$-complete, $p(B_c)\geq4$, $\emph{Sing}(B_+)=\ES$, and $\xi$ be a divisor on $C_+$ with zero Chern class in $H^2(C_+,\Z)$. Then $\xi$ extends (uniquely) to a divisor on $B_+$.
\ec
\demo
Use diagram~(\ref{last}).
\enddemo

\bt\label{ansets}
Assume that $H^2(C_+,\Q)=\{0\}$. If $\emph{Sing}(B_+)=\ES$, $B_c$ is $1$-complete and $p(B_c)\geq4$, then all analytic sets of codimension $1$ on $C_+$ extend to analytic sets on $B_+$.
\et
\demo
Let $A$ be an analytic set of codimension $1$ on $C_+$. Since $B_+$ is a Stein manifold, $C_+$ is locally factorial, and so there exists a divisor $\xi$ on $C_+$ with support $A$. Since 
$H^2(C_+,\Q)=\{0\}$, there exists $n\in\N$ such that $n c_2(\xi)=0\in H^2(C_+,\Z)$. Hence $n\xi$ has zero Chern class in $H^2(C_+,\Z)$, and so, by Corollary~\ref{corxi} $n\xi$ can be extended to a divisor $\widetilde{n\xi}$ on $B_+$. The support of $\widetilde{n\xi}$ is an analytic set $\widetilde A$ which extends to $B_+$ the support $A$ of $n\xi$.
\enddemo

In Theorem~\ref{divis} the condition $H^2(C_+,\Z)=\{0\}$ can be relaxed and replaced by the weaker one: the restriction map $H^2(B_+,\Z)\to H^2(C_+,\Z)$ is surjective. We need the following
\begin{Lemma}\label{deltaI}
$\delta$ is injective.
\end{Lemma}
\demo
First we prove lemma for $C_+$ closed. Let $\xi\in H^1(\oli B_+,\mathcal M^*)$ be such that $\xi_{|\oli C_+}=0$. Consider the set
$$A=\{\eta\in[0,\varepsilon]\ :\  \xi_{|\oli B_+\smi\oli B_{\eta}}=0\}.$$
If we prove that $0\in A$, we are done, because $0=\xi_{|\oli B_+\smi \oli B_0}=\xi_{|\oli B_+}=\xi$. Obviously $\eta_0\in A$ implies $\forall\eta>\eta_0$, $\eta\in A$.

$A\neq\ES$. Since $C_+=B_+\smi \oli B_\e$ and $\xi_{|\oli C_+}=0$, $\varepsilon\in A$.

$A$ is closed. If $\eta_n\in A$, for all $n$, and $\eta_n\searrow\eta_\infty$, $\oli B_+\smi \oli B_{\eta_\infty}=\cup_n (\oli B_+\smi \oli B_{\eta_n})$, hence $\xi_{|\oli B_+\smi \oli B_{\eta_n}}=0$ for all $n$ implies $\xi_{|\oli B_+\smi \oli B_{\eta_\infty}}=0$, i.e.\ $\eta_\infty\in A$.

$A$ is open. Suppose $0<\eta_0\in A$. We denote $C_{\eta_0}=\oli B_+\smi \oli B_{\eta_0}$. Let $\mathcal A$ be the family of open covering $\{U_i\}_{i\in I}$ of $\oli B_+$ such that:
\begin{itemize}
\item[$\alpha$)] $U_i$ is isomorphically equivalent to an holomorphy domain in $\C^n$;
\item[$\beta$)] If $U_i\cap b B_{\eta_0}\neq\ES$, the restriction homomorphism
$$
H^0(U_i,\mathcal O)\rightarrow H^0(U_i\cap C_{\eta_0},\mathcal O)$$
is bijective;
\item[$\gamma$)] $U_i\cap U_j$ is simply connected.
\end{itemize}
$\mathcal A$ is not empty and it is cofinal in the set of open coverings of $\oli B_+$ \cite[Lemma 2, p.\ 222]{AG}. Let $\mathcal U=\{U_i\}_{i\in I}\in \mathcal A$, and $\{f_{ij}\}\in Z^1(\mathcal U,\mathcal M^*)$ be a representative of $\xi$. Let $W_i=U_i\cap C_{\eta_0}$. Since $\eta_0\in A$, if $W_i\cap W_j\neq\ES$, $f_{ij|W_i\cap W_j}=f_i f_j^{-1}$ ($f_\nu\in H^0(W_\nu,\mathcal M^*)$). By $\alpha$), $f_\nu=p_\nu q_\nu^{-1}$, $p_\nu, q_\nu\in H^0(W_\nu,\mathcal O)$. By $\beta$), both $p_\nu$ and $q_\nu$ can be holomorphically extended on $U_\nu$, with $\widetilde p_\nu$ and $\widetilde q_\nu$. Hence we have $f_{ij}=\widetilde p_i \widetilde q_i^{-1}(\widetilde p_j \widetilde q_j^{-1})^{-1}$ on $U_i\cap U_j$ (which is simply connected, so that there is no polidromy). So $\xi=0$ in an open neighborhood $U$ of $C_{\eta_0}$ and, by compactness, there exists $\epsilon'>0$ such that $C_{\eta_0-\epsilon'}\subset U$. So $\eta_0-\epsilon'\in A$ and consequently $A$ is open.

Thus $A=[0,\varepsilon]$, and the lemma is proved if $C_+$ is closed.

If $C_+$ is open, we consider $C_+$ as a union of the closed semi $1$-coronae
$$
\oli C_\epsilon=\oli{B_{\varepsilon+\epsilon',c}\cap\{h>\epsilon'\}}\subset C_+.
$$
Let $\xi\in H^1(B_+,\mathcal M^*)$ be such that $\xi_{|C_+}=0$. Then $\xi_{|\oli C_\epsilon'}=0$, for all $\epsilon'>0$. Consequently from what we have already proved $\xi_{|\oli B_\epsilon'}=0$, where $\oli B_\epsilon=\oli{B_+\cap\{h>\epsilon'\}}$. Since $\cup_\epsilon' \oli B_\epsilon'=B_+$, $\xi=0$ and the lemma is proved.
\enddemo

Lemma~\ref{alphaS}, Lemma~\ref{gammaS} and Lemma~\ref{deltaI} lead to the following generalization of Theorem~\ref{divis}:
\bt\label{divis2}
Assume that the restriction $H^2(B_+,\Z)\to H^2(C_+,\Z)$ is surjective. If $\emph{Sing}(B_+)=\ES$, $B_c$ is $1$-complete and $p(B_c)\geq4$, then all divisors on $C_+$ extend to divisors on $B_+$.
\et

\end{document}